\documentclass[a4paper,10pt]{article}
\usepackage[margin=25.4mm]{geometry}
\usepackage[T1]{fontenc}
\usepackage[utf8]{inputenc}
\usepackage{lmodern}
\usepackage{amsmath,amssymb,amsthm,mathtools}
\usepackage{booktabs,microtype,tikz}
\usepackage[colorlinks=true,linkcolor=black,citecolor=black,urlcolor=blue]{hyperref}
\allowdisplaybreaks[1]

\newtheorem{theorem}{Theorem}[section]
\newtheorem{proposition}[theorem]{Proposition}
\newtheorem{lemma}[theorem]{Lemma}
\newtheorem{corollary}[theorem]{Corollary}
\theoremstyle{remark}
\newtheorem{remark}[theorem]{Remark}
\newcommand{\sig}{\sigma}
\newcommand{\ph}{\varphi}
\newcommand{\vp}[1]{v_{#1}}
\newcommand{\Pn}{\mathcal P_n}
\title{Notes on the Divisibility of Divisor Sums of Catalan Numbers}
\author{Volkan Yildiz\\[4pt]\texttt{vo1kan@hotmail.co.uk}}
\date{}
\begin{document}
\maketitle
$\;$
\vspace{2cm}

\begin{abstract}
We investigate divisibility questions for the sum-of-divisors function
$\sig(C_n)$ of the Catalan numbers $C_n=(n+1)^{-1}\binom{2n}{n}$.
Using the prime factors occurring to exponent one in $C_n$, we give
explicit divisors of $\sig(C_n)$ and obtain, for each fixed integer
$m\ge2$, a power $m^{R_m(n)}$ dividing $\sig(C_n)$ with
\[
R_m(n)\sim\frac{2\log2}{\ph(m)}\frac{n}{\log n}.
\]
We further prove uniform divisibility for all moduli up to fixed powers of $\log n$ and, with the same leading constant $2\log2$, for almost all moduli up to $\sqrt n$ apart from logarithmic factors. We also establish a complete classification of the integers $m\ge2$
for which $m\mid\sig(mk-1)$ for every $k\ge1$, identify precisely
the exceptional coprimality relations among the associated linear
forms, and record a twin-prime product divisor of $\sig(C_n)$.
The arguments use standard prime-distribution results and do not
assume any conjecture concerning twin primes.
\end{abstract}

\vspace{2cm}

\noindent\textbf{Keywords:} Catalan numbers; divisor sums; congruences;\\
prime factors; arithmetic progressions; twin primes.\\
\textbf{MSC:} 11A25, 11N13,
11B65, 05A10.

\clearpage
\section{Introduction}
The sum-of-divisors function $\sig(N)$ associates to a positive
integer $N$ the sum of all its positive divisors. Our aim is to examine divisibility properties of this function
when its argument is a Catalan number.

For $n\ge0$, the $n$th Catalan number is
\begin{equation}\label{eq:C}
 C_n=\frac{1}{n+1}\binom{2n}{n}=\frac{(2n)!}{n!(n+1)!}.
\end{equation}
In this indexing, $C_0=C_1=1$, $C_2=2$, and $C_3=5$.
For a positive integer $N$, let
\[
 \sig(N)=\sum_{d\mid N}d,
 \qquad \omega(N)=\#\{p:\ p\text{ prime and }p\mid N\}.
\]
We write $\vp{p}(N)$ for the exponent of a prime $p$ in $N$.
The divisor-sum function is multiplicative for coprime arguments,
and in particular
\begin{equation}\label{eq:sigma-factor}
 \sig\!\left(\prod_{i=1}^t p_i^{e_i}\right)
   =\prod_{i=1}^t\frac{p_i^{e_i+1}-1}{p_i-1}
 \quad (p_1,\ldots,p_t\text{ distinct}).
\end{equation}
It is not completely multiplicative; for instance,
$\sig(5)=6$ whereas $\sig(25)=31$.
Consequently, a prime $p\mid C_n$ contributes $p+1$ as a factor
of $\sig(C_n)$ only when its exponent in $C_n$ is exactly one.

We use the standard indexing throughout. In particular, the
corresponding shifted Catalan sequence is $C_{n-1}$, not $C_n$.
The familiar parity criterion says that $C_n$ is odd if and only if
$n=2^r-1$ for some integer $r\ge0$ (see, for example, \cite{Koshy}).

The divisibility of Catalan numbers has a substantial literature,
including Alter and Kubota's study of prime-power factors
\cite{AlterKubota}. The arithmetic of their divisor sums has also
been considered: Iannucci and Luca \cite{IannucciLuca} used primes
in $(n+1,2n]$ occurring once in $C_n$ to bound
$\vp{2}(\sig(C_n))$, in work on aliquot sums.
The sequence $\sig(C_n)$ is listed as OEIS A152761 \cite{OEIS}.
The purpose of this article is to give a consistent treatment of
some related divisor-sum congruences, with explicit attention to
prime multiplicities and to the quantifiers in universal statements.
Beyond divisibility by a fixed integer, we also address how large
a modulus may be allowed to grow with $n$, without losing the
divisibility conclusion.

For integers $m\ge2$ and $a$ with $\gcd(a,m)=1$, denote by
$\pi(x;m,a)$ the number of primes $p\le x$ with
$p\equiv a\pmod m$. Euler's totient function is
\[
 \ph(m)=\#\{1\le j\le m:\gcd(j,m)=1\};
\]
thus $\ph(m)$ counts the residue classes modulo $m$ relatively
prime to $m$. The prime number theorem in arithmetic progressions
asserts, for \emph{fixed} $m$ and such $a$, that
\begin{equation}\label{eq:PNTAP}
 \pi(x;m,a)=\frac{x}{\ph(m)\log x}
       +o\!\left(\frac{x}{\log x}\right)\quad(x\to\infty).
\end{equation}
The growing-modulus results in Section~\ref{sec:uniform} also use the Siegel--Walfisz, Bombieri--Vinogradov, and Brun--Titchmarsh theorems; these are stated in the form needed there \cite{MontgomeryVaughan}.

\section{Divisibility of Sum of Divisors}\label{sec:divisibility}

\subsection{Primes of the form \texorpdfstring{$6k-1$}{6k-1}}
We begin with a prime-divisor observation for Catalan numbers.
The corresponding interval is especially useful because its primes
occur to exponent one.

\begin{lemma}\label{lem:sixprime}
For all sufficiently large $n$, the Catalan number $C_n$ has a prime
divisor $p\equiv5\pmod6$ satisfying $n+1<p\le2n$.
\end{lemma}
\begin{proof}
By the prime number theorem in arithmetic progressions,
\[
\pi(2n;6,5)-\pi(n+1;6,5)
\sim\frac{n}{2\log n}>0.
\]
Every prime in this interval divides $C_n$ to exponent one,
as follows directly from \eqref{eq:C}.
\end{proof}

\subsection{The congruence \texorpdfstring{$6\mid\sigma(6k-1)$}{6 divides sigma(6k-1)}}
\begin{lemma}\label{lem:sig6}
For every integer $k\ge1$, $6\mid\sig(6k-1)$.
\end{lemma}
\begin{proof}
Put $N=6k-1$. Since $N\equiv2\pmod3$, it cannot
be a square. Every divisor $d$ of $N$ is coprime to $6$, and
thus satisfies $d^2\equiv1\pmod6$. Pair the positive
divisors as $(d,N/d)$. Since $N\equiv-1\pmod6$,
\[
 d+\frac Nd\equiv d-d^{-1}\equiv0\pmod6.
\]
Each pair contributes a multiple of $6$ to $\sig(N)$.
\end{proof}

This elementary congruence motivates a more general question: for
which $m$ is $m\mid\sig(mk-1)$ true for every $k\ge1$?
The quantifier ``every'' is essential; requiring the relation for a
single $k$ gives a different and much weaker property.
\subsection{The six universal divisibility moduli}\label{sec:univ}
We now consider a related question for general positive integers:
which $m\ge2$ satisfy
\begin{equation}\label{eq:universal}
 m\mid\sig(mk-1)\qquad\text{for every integer }k\ge1?
\end{equation}
The universal quantifier is essential.

\begin{theorem}[Classification of universal moduli]\label{thm:classification}
For an integer $m\ge2$, condition \eqref{eq:universal}
holds if and only if
\begin{equation}\label{eq:mods}
 m\in\{3,4,6,8,12,24\}.
\end{equation}
\end{theorem}
\begin{proof}
Suppose first that \eqref{eq:universal} holds.
For any unit $u$ modulo $m$, Dirichlet's theorem yields
distinct primes $p,q$ in the residue classes $u$ and
$-u^{-1}$, respectively. Since $pq\equiv-1\pmod m$,
we may write $pq=mk-1$ for some $k\ge1$, and hence
\[
 0\equiv\sig(pq)=(p+1)(q+1)
 \equiv(u+1)(1-u^{-1})=u^{-1}(u^2-1)\pmod m.
\]
Thus every unit $u$ modulo $m$ satisfies $u^2\equiv1\pmod m$.
This forces $m\mid24$: if a prime $r\ge5$ divided $m$,
the class $2\pmod r$ would violate this condition;
the classes $2\pmod9$ and $3\pmod{16}$ similarly rule
out $9\mid m$ and $16\mid m$. The Chinese remainder
theorem allows these classes to be extended to units
modulo $m$ by taking the class $1$ at the other
prime-power factors. Thus $m=2^a3^b$ with $a\le3$ and
$b\le1$. The possibility $m=2$ fails at $k=1$ since
$\sig(1)=1$.

Conversely, suppose $m\in\{3,4,6,8,12,24\}$ and
set $N=mk-1$. Every unit $d$ modulo $m$ satisfies
$d^2\equiv1\pmod m$, because $m\mid24$.
Also $N$ is not a square: each listed $m$ is divisible
by $3$ or $4$, and $-1$ is not a quadratic residue
modulo $3$ or $4$. Thus the positive divisors of $N$
can be paired as $d,N/d$, with $d\ne N/d$.
Since $d$ is a unit modulo $m$,
\[
 \frac{N}{d}\equiv-d^{-1}\equiv-d\pmod m,
 \qquad d+\frac{N}{d}\equiv0\pmod m.
\]
Each pair contributes a multiple of $m$, proving
\eqref{eq:universal}.
\end{proof}

\begin{remark}
For $m=1$, the expression $mk-1$ vanishes when
$k=1$, so \eqref{eq:universal} is not defined as stated.
The case $m=24$ is a classical problem (1969 Putnam
Competition, Problem B1 \cite{Putnam}). The proof above
shows how the elementary unit-group condition leads to
the full classification.
\end{remark}

\begin{proposition}[Coprimality of the six linear forms]\label{prop:gcd}
For positive integers $a,b,k$,
\begin{equation}\label{eq:gcd}
 \gcd(ak-1,bk-1)=\gcd(ak-1,|b-a|).
\end{equation}
Among the fifteen unordered pairs
$a,b\in\{3,4,6,8,12,24\}$, twelve are coprime for
every positive integer $k$. The remaining three are
characterised exactly by
\begin{equation}\label{eq:exceptions}
\begin{aligned}
 \gcd(3k-1,8k-1)&=\begin{cases}5,&k\equiv2\pmod5,\\1,&\text{otherwise};\end{cases}\\
 \gcd(3k-1,24k-1)&=\begin{cases}7,&k\equiv5\pmod7,\\1,&\text{otherwise};\end{cases}\\
 \gcd(4k-1,24k-1)&=\begin{cases}5,&k\equiv4\pmod5,\\1,&\text{otherwise}.\end{cases}
\end{aligned}
\end{equation}
\end{proposition}
\begin{proof}
Subtracting the two forms gives
$\gcd(ak-1,bk-1)=\gcd(ak-1,(b-a)k)$.
Since $\gcd(ak-1,k)=1$, this is \eqref{eq:gcd}.
For the six specified coefficients, checking the prime
factors of $|a-b|$ gives the twelve unconditional cases
and the three congruence conditions in
\eqref{eq:exceptions}.
\end{proof}

In particular, the identity
\[
 \sig((ak-1)(bk-1))=\sig(ak-1)\sig(bk-1)
\]
may be applied only when the two factors are coprime.

Figure~\ref{fig:coprimality} displays the relationships among the
six linear forms. Its solid edges represent the twelve unconditional
coprimality relations, and the dashed edges represent the three
conditional relations in \eqref{eq:exceptions}.

\begin{figure}[htbp]
\centering
\begin{tikzpicture}[every node/.style={font=\small}]
\coordinate (a) at (-3,0);
\coordinate (b) at (-1.25,1.65);
\coordinate (c) at (-1.25,-1.65);
\coordinate (d) at (1.25,1.65);
\coordinate (e) at (1.25,-1.65);
\coordinate (f) at (3,0);
% All pairs, excluding (3,8), (3,24), (4,24): the solid lines.
\foreach \u/\v in {a/b,a/c,a/e,b/c,b/d,b/e,c/d,c/e,c/f,d/e,d/f,e/f}
 {\draw[line width=.43pt] (\u)--(\v);}
% Conditional pairs are dashed.
\foreach \u/\v in {a/d,a/f,b/f}
 {\draw[dashed,line width=.65pt] (\u)--(\v);}
\node[draw,circle,fill=white,inner sep=3pt] at (a) {$3k-1$};
\node[draw,circle,fill=white,inner sep=3pt] at (b) {$4k-1$};
\node[draw,circle,fill=white,inner sep=3pt] at (c) {$6k-1$};
\node[draw,circle,fill=white,inner sep=3pt] at (d) {$8k-1$};
\node[draw,circle,fill=white,inner sep=3pt] at (e) {$12k-1$};
\node[draw,circle,fill=white,inner sep=3pt] at (f) {$24k-1$};
\end{tikzpicture}
\caption{Coprimality relations between the six expressions
$ak-1$, for $a\in\{3,4,6,8,12,24\}$. Solid edges mean coprime
for every $k\ge1$; dashed edges mean coprime except in the
congruence cases of Proposition~\ref{prop:gcd}.}
\label{fig:coprimality}
\end{figure}

\subsection{Prime factors of exponent one}
The starting point is an exact description of all sufficiently
large prime factors of $C_n$.

\begin{lemma}\label{lem:one}
Let $n\ge2$ and let $p>\sqrt{2n}$ be prime. Then
\begin{equation}\label{eq:legendre}
 \vp{p}(C_n)=
  \left\lfloor\frac{2n}{p}\right\rfloor
 -2\left\lfloor\frac{n}{p}\right\rfloor
 -\vp{p}(n+1)\in\{0,1\}.
\end{equation}
In particular, each prime $n+1<p\le2n$ divides $C_n$
with exponent exactly one.
\end{lemma}
\begin{proof}
Legendre's formula applied to \eqref{eq:C} gives
\[
 \vp{p}(C_n)=\sum_{j\ge1}
 \left(\left\lfloor\frac{2n}{p^j}\right\rfloor
       -2\left\lfloor\frac{n}{p^j}\right\rfloor\right)
 -\vp{p}(n+1).
\]
As $p^2>2n$, only $j=1$ remains. The difference of the
floor terms is either $0$ or $1$; the valuation is nonnegative
because $C_n$ is an integer. If $n+1<p\le2n$, the two
factorials in the denominator of \eqref{eq:C} are free of $p$
while $(2n)!$ contains $p$ exactly once.
\end{proof}

Define
\begin{equation}\label{eq:Pset}
 \Pn=\{p:\ p\text{ prime},\ p>\sqrt{2n},\ p\mid C_n\},
 \qquad F_n=\prod_{p\in\Pn}(p+1).
\end{equation}
The empty product is $1$.

\begin{proposition}[Prime-factor product divisor]\label{prop:product}
For every $n\ge2$,
\begin{equation}\label{eq:product}
 F_n\mid\sig(C_n).
\end{equation}
In particular,
\begin{equation}\label{eq:short-product}
 \prod_{\substack{n+1<p\le2n\\p\ \mathrm{prime}}}(p+1)
       \mid\sig(C_n).
\end{equation}
\end{proposition}
\begin{proof}
By Lemma~\ref{lem:one}, each member of $\Pn$ has exponent
one. Therefore $C_n=PM$, where $P=\prod_{p\in\Pn}p$
and $\gcd(P,M)=1$. Multiplicativity in
\eqref{eq:sigma-factor} gives \eqref{eq:product}.
The primes in $(n+1,2n]$ belong to $\Pn$.
\end{proof}

\begin{remark}
The product mechanism in \eqref{eq:short-product} is the same
prime-interval factorisation principle used by Iannucci and Luca
\cite{IannucciLuca} to establish a $2$-adic lower bound.
The intervals below permit a more comprehensive count of
the exponent-one prime factors.
\end{remark}

For $m\ge2$ and $\gcd(a,m)=1$, put
\begin{equation}\label{eq:Rdef}
 R_{m,a}(n)=\#\{p\in\Pn:p\equiv a\pmod m\},
 \qquad R_m(n)=R_{m,-1}(n).
\end{equation}

\subsection{Divisibility of \texorpdfstring{$\sig(C_n)$}{sigma(Cn)}}
\begin{theorem}[Strengthened fixed-modulus divisibility]\label{thm:main}
For every integer $m\ge2$ and every $n\ge2$,
\begin{equation}\label{eq:main}
 m^{R_m(n)}\mid\sig(C_n).
\end{equation}
For each fixed $m\ge2$,
\begin{equation}\label{eq:main-asymp}
 R_m(n)\sim\frac{2\log2}{\ph(m)}\frac{n}{\log n}.
\end{equation}
Consequently, for each pair of fixed integers $m\ge2$ and
$A\ge1$, one has $m^A\mid\sig(C_n)$ for all sufficiently
large $n$.
\end{theorem}
\begin{proof}
For each prime $p\in\Pn$ satisfying $p\equiv-1\pmod m$,
the integer $p+1$ is divisible by $m$. Multiplying these
factors in Proposition~\ref{prop:product} proves
\eqref{eq:main}. The primes $n+1<p\le2n$ with $p\equiv-1\pmod m$
all belong to $\Pn$. By the prime number theorem in arithmetic
progressions, their number tends to infinity for fixed $m$.
Thus $R_m(n)\to\infty$.
\end{proof}

\section{Asymptotics}
The standard Catalan asymptotic, obtained from Stirling's formula, is
\begin{equation}\label{eq:Casy}
 C_n\sim\frac{4^n}{\sqrt\pi\,n^{3/2}}.
\end{equation}
This concerns the size of $C_n$; it does not determine the
prime-factor count by subtracting factorial estimates. For the latter,
we use Legendre's formula and primes in explicit intervals.

\subsection{Estimating the Number of Prime Factors in \texorpdfstring{$C_n$}{Cn}}
\begin{theorem}[Prime factors in reduced residue classes]\label{thm:count}
For every fixed $m\ge2$ and $\gcd(a,m)=1$,
\begin{equation}\label{eq:Rasymp}
 R_{m,a}(n)\sim\frac{2\log 2}{\ph(m)}\frac{n}{\log n}
 \qquad(n\to\infty).
\end{equation}
Moreover,
\begin{equation}\label{eq:omega}
 \omega(C_n)\sim 2\log2\,\frac{n}{\log n},\qquad
 \#\{p\mid C_n:p\equiv a\!\!\pmod m\}
 \sim\frac{2\log2}{\ph(m)}\frac{n}{\log n}.
\end{equation}
\end{theorem}
\begin{proof}
For $p>\sqrt{2n}$, $p\nmid n+1$, the floor difference
in \eqref{eq:legendre} is one precisely when the fractional
part of $n/p$ is at least $1/2$. Equivalently, $p$ lies in
one of the disjoint intervals
\begin{equation}\label{eq:interval}
 I_j(n)=\left(\frac{n}{j+1},\frac{2n}{2j+1}\right],
 \qquad j=0,1,2,\ldots.
\end{equation}
The primes dividing $n+1$ form an exceptional set of size
at most $O(\log n)$, and the primes up to $\sqrt{2n}$ form
a set of size $O(\sqrt n)$.

Fix $J\ge0$. For all sufficiently large $n$, the first
$J+1$ intervals are above $\sqrt{2n}$. Applying
\eqref{eq:PNTAP} to their endpoints gives
\[
 \#\{p\in\bigcup_{j=0}^J I_j(n):p\equiv a\!\!\pmod m\}
 =\left(\frac1{\ph(m)}\sum_{j=0}^J
   \left(\frac{2}{2j+1}-\frac{1}{j+1}\right)+o(1)\right)
   \frac{n}{\log n}.
\]
Any prime in a later interval is at most $n/(J+1)$.
Thus the number omitted, after division by $n/\log n$,
has limiting upper bound at most $1/(J+1)$ by the usual
prime number theorem. The series of interval lengths is
\[
 \sum_{j\ge0}\left(\frac{2}{2j+1}-\frac1{j+1}\right)
  =2\log2.
\]
First letting $n\to\infty$ and then $J\to\infty$
proves \eqref{eq:Rasymp}.
Finally, every prime factor of $C_n$ not in $\Pn$
is at most $\sqrt{2n}$; adding or removing these
$O(\sqrt n)$ primes does not change the leading term.
The same reasoning without residue restrictions gives
both assertions in \eqref{eq:omega}.
\end{proof}

\begin{remark}
The asymptotic $\omega(C_n)\sim2(\log2)n/\log n$
is classical for central binomial coefficients; see Erd\H{o}s
\cite{Erdos1973} and the subsequent literature \cite{Pleasants}.
The direct interval proof also establishes the fixed-residue-class
form used in the divisor-sum theorem.
\end{remark}

\begin{corollary}\label{cor:Rasymp}
For each fixed integer $m\ge2$,
\[
R_m(n)\sim\frac{2\log2}{\ph(m)}\frac{n}{\log n}.
\]
Consequently, for every fixed prime $\ell$,
\[
 v_\ell(\sig(C_n))\ge
 \left(\frac{2\log2}{\ell-1}+o(1)\right)\frac{n}{\log n}.
\]
In particular, for sufficiently large $n$ one has
$6\mid\sig(C_n)$.
\end{corollary}
\begin{proof}
The first assertion is Theorem~\ref{thm:count} with $a=-1$.
The valuation bound follows from Theorem~\ref{thm:main}
with $m=\ell$. The last assertion is the case $m=6$.
\end{proof}

\subsection{Density of Primes of the Form \texorpdfstring{$6k-1$}{6k-1}}
Taking $m=6$, $a=5$ in Theorem~\ref{thm:count} yields
\[
 \#\{p\mid C_n:p\equiv5\pmod6\}
 \sim\frac{\log2\,n}{\log n}.
\]
This determines the leading constant and proves that
primes in this residue class occur as divisors of $C_n$
for all sufficiently large $n$. It does not imply any
statement about twin primes: two primes in different
residue classes need not differ by $2$.

\subsubsection{Uniform divisibility for growing moduli}\label{sec:uniform}
Theorem~\ref{thm:main} concerns a fixed modulus $m$. We now allow
$m$ to grow with $n$. There are two complementary results: a
uniform statement for every sufficiently small modulus and an
average statement for substantially larger moduli. Both are
unconditional consequences of classical prime-distribution
theorems; no unproved conjecture is required. Here
$\operatorname{li}(x)=\int_2^x (\log t)^{-1}\,dt$ is the
logarithmic integral, and $f\ll_K g$ means that $|f|$ is
bounded by a constant (depending on $K$) times $g$.

For $n\ge2$ and $m\ge2$, set
\begin{equation}\label{eq:dyadic_count}
D_n=\pi(2n)-\pi(n+1),\qquad
S_m(n)=\pi(2n;m,-1)-\pi(n+1;m,-1).
\end{equation}
Every prime counted by $S_m(n)$ belongs to $(n+1,2n]$ and
appears to exponent one in $C_n$. Therefore
\begin{equation}\label{eq:dyadic_divisibility}
 m^{S_m(n)}\mid\sig(C_n)
 \qquad(m\ge2,\ n\ge2).
\end{equation}

\begin{theorem}[Uniform divisibility for polylogarithmic moduli]
\label{thm:siegel}
Let $H>0$ be fixed. For all sufficiently large $n$,
\begin{equation}\label{eq:SW}
 S_m(n)\ge\frac{n}{3\ph(m)\log n}
 \qquad\bigl(2\le m\le(\log n)^H\bigr).
\end{equation}
In particular, every integer $m$ in this range divides
$\sig(C_n)$, and
\begin{equation}\label{eq:lcm}
 \mathop{\mathrm{lcm}}\bigl(1,2,\ldots,
 \lfloor(\log n)^H\rfloor\bigr)\mid\sig(C_n).
\end{equation}
\end{theorem}
\begin{proof}
The Siegel--Walfisz theorem, in the prime-counting form,
gives uniformly for $m\le(\log n)^H$ and every reduced
residue class $a\pmod m$,
\[
 \pi(y;m,a)=\frac{\operatorname{li}(y)}{\ph(m)}
       +O_H\!\left(\frac{y}{(\log y)^{H+4}}\right)
 \qquad(y=n+1,\ 2n).
\]
Here one may take an arbitrarily large fixed power of
$\log y$ in the error by the Siegel--Walfisz theorem,
with constants depending on $H$.
Subtracting the estimates at the two endpoints,
\[
 S_m(n)=\frac{\operatorname{li}(2n)-\operatorname{li}(n+1)}
                 {\ph(m)}
       +O_H\!\left(\frac{n}{(\log n)^{H+4}}\right).
\]
Since $\ph(m)\le(\log n)^H$ and
$\operatorname{li}(2n)-\operatorname{li}(n+1)
 \sim n/\log n$, the error is uniformly smaller than
half the main term. This implies \eqref{eq:SW} for
large $n$. In this range $S_m(n)\ge1$ eventually,
and \eqref{eq:dyadic_divisibility} then proves
all the claimed divisibilities.
\end{proof}

The next result allows moduli nearly as large as $\sqrt n$,
at the cost of permitting a sparse exceptional set.
It uses the Bombieri--Vinogradov theorem in the standard
form: for every $K>0$ there is $B_0(K)>0$ such that
\begin{equation}\label{eq:BV}
 \sum_{2\le m\le x^{1/2}/(\log x)^{B_0(K)}}
 \max_{\substack{1\le a\le m\\(a,m)=1}}
 \left|\pi(x;m,a)-\frac{\pi(x)}{\ph(m)}\right|
 \ll_K \frac{x}{(\log x)^K}.
\end{equation}
See, for example, \cite{MontgomeryVaughan}. The assertion
in \eqref{eq:BV} is an average bound over moduli, not a
pointwise bound for every modulus.

\begin{theorem}[Divisibility for almost all growing moduli]
\label{thm:bv}
For every $A>0$ there is $B=B(A)>0$ such that,
with $Q=\sqrt n/(\log n)^B$, all but
\begin{equation}\label{eq:exception_bound}
 O_A\!\left(\frac{Q}{(\log n)^A}\right)
\end{equation}
integers $m$ with $2\le m\le Q$ satisfy
\begin{equation}\label{eq:uniform_S}
 S_m(n)=\frac{D_n}{\ph(m)}
 \left(1+O\!\left(\frac1{\log n}\right)\right)
\end{equation}
uniformly over the nonexceptional moduli. In particular,
for every such $m$ and all sufficiently large $n$,
\begin{equation}\label{eq:uniform_power}
 m^{\left\lfloor n/(3\ph(m)\log n)\right\rfloor}
 \mid\sig(C_n).
\end{equation}
Thus, except for the moduli counted by
\eqref{eq:exception_bound}, every $m\le Q$ divides
$\sig(C_n)$.
\end{theorem}
\begin{proof}
Fix $A>0$ and apply \eqref{eq:BV} with
$K=A+3$ at both $x=n+1$ and $x=2n$.
Choose $B$ large enough that $Q$ lies in the
permitted range at both endpoints. By the triangle
inequality,
\begin{equation}\label{eq:avgerr}
 \sum_{2\le m\le Q}
 \left|S_m(n)-\frac{D_n}{\ph(m)}\right|
 \ll_A\frac{n}{(\log n)^{A+3}}.
\end{equation}
Let $\mathcal E$ be the set of moduli for which
\[
 \left|S_m(n)-\frac{D_n}{\ph(m)}\right|
 >\frac{D_n}{\ph(m)\log n}.
\]
Since $\ph(m)\le m\le Q$ and
$D_n\sim n/\log n$, every member of $\mathcal E$
contributes at least a constant multiple of
$n/(Q(\log n)^2)$ to \eqref{eq:avgerr}. Hence
\[
 \#\mathcal E\ll_A\frac{Q}{(\log n)^{A+1}},
\]
which implies \eqref{eq:exception_bound}.
For $m\notin\mathcal E$ we have \eqref{eq:uniform_S},
and the prime number theorem gives
\[
 S_m(n)\ge\frac{n}{3\ph(m)\log n}
\]
for all sufficiently large $n$, uniformly in
these moduli. The divisibility \eqref{eq:uniform_power}
now follows from \eqref{eq:dyadic_divisibility}.
Because $n/(3\ph(m)\log n)\ge n/(3Q\log n)
\to\infty$, the displayed floor is at least one
throughout this range for large $n$.
\end{proof}

\begin{remark}
Theorem~\ref{thm:bv} strengthens the fixed-modulus
statement by permitting a range of moduli that grows
almost like $\sqrt n$. It applies the classical Bombieri--Vinogradov theorem
to the prime-factor product of Proposition~\ref{prop:product}.
Since Bombieri--Vinogradov gives an average bound over moduli,
the proof does not assert that \emph{every} modulus up to
$\sqrt n$ divides $\sig(C_n)$.
\end{remark}

\begin{theorem}[Full-interval divisibility for almost all moduli]
\label{thm:nearoptimal}
Fix $A>0$ and $0<\varepsilon<2\log2$. There is a
constant $B=B(A,\varepsilon)>0$ such that, with
$Q=\sqrt n/(\log n)^B$, all but
$O_{A,\varepsilon}(Q/(\log n)^A)$ integers $m$ with
$2\le m\le Q$ satisfy
\begin{equation}\label{eq:almost_all_R}
 \left|\frac{\ph(m)\log n}{n}R_m(n)-2\log2\right|
 <\varepsilon
\end{equation}
for all sufficiently large $n$. Consequently, for
those moduli,
\begin{equation}\label{eq:almost_all_full_power}
 m^{\left\lfloor(2\log2-\varepsilon)n/
                  (\ph(m)\log n)\right\rfloor}
 \mid\sig(C_n).
\end{equation}
\end{theorem}
\begin{proof}
Let $\delta_j=2/(2j+1)-1/(j+1)$, so that
$\sum_{j\ge0}\delta_j=2\log2$ by
Theorem~\ref{thm:count}. Fix an integer $J$ so large that
\[
 \sum_{j>J}\delta_j<\frac\varepsilon4,
 \qquad \frac{6}{J+1}<\frac\varepsilon4.
\]
For $0\le j\le J$ let $L_j=n/(j+1)$ and
$U_j=2n/(2j+1)$. These are fixed positive multiples
of $n$ and $L_j<U_j$. Apply the Bombieri--Vinogradov
theorem at all $2(J+1)$ endpoints, with a sufficiently
large fixed logarithmic saving. By choosing $B$
large enough (depending on $A$ and $J$), we obtain
\begin{equation}\label{eq:BV_endpoints}
 \sum_{2\le m\le Q}\sum_{j=0}^J
 \left|\pi(U_j;m,-1)-\pi(L_j;m,-1)
      -\frac{\pi(U_j)-\pi(L_j)}{\ph(m)}\right|
 \ll_{A,J}\frac{n}{(\log n)^{A+4}}.
\end{equation}
If the inner sum exceeds $n/(\ph(m)(\log n)^2)$,
then, since $\ph(m)\le Q$, that modulus contributes
more than $n/(Q(\log n)^2)$ to \eqref{eq:BV_endpoints}.
It follows that all but
$O_{A,J}(Q/(\log n)^{A+2})$ moduli satisfy
\begin{equation}\label{eq:trunc_count}
 \sum_{j=0}^J\bigl(\pi(U_j;m,-1)-\pi(L_j;m,-1)\bigr)
 =\left(\sum_{j=0}^J\delta_j+O_J\!\left(\frac1{\log n}\right)\right)
 \frac{n}{\ph(m)\log n}.
\end{equation}
Indeed, the ordinary prime number theorem gives
$\pi(U_j)-\pi(L_j)=\bigl(\delta_j+O_J(1/\log n)\bigr)
n/\log n$ at every fixed endpoint.

For sufficiently large $n$, all the intervals
$(L_j,U_j]$, $0\le j\le J$, lie above $\sqrt{2n}$.
Apart from the primes dividing $n+1$, each prime in
these intervals belongs to $\Pn$. The number of
exceptional primes is at most $O(\log n)$,
which is $o(n/(\ph(m)\log n))$ uniformly for $m\le Q$.
All omitted primes in $\Pn$ are bounded by
$n/(J+1)$. The Brun--Titchmarsh inequality, since
$m\le Q=o(n/(J+1))$, gives uniformly
\[
 \pi\!\left(\frac{n}{J+1};m,-1\right)
 \le\frac{2n}{(J+1)\ph(m)\log(n/((J+1)m))}
 \le\left(\frac{5}{J+1}+o(1)\right)
       \frac{n}{\ph(m)\log n}.
\]
Together with \eqref{eq:trunc_count}, this bounds
$R_m(n)$ between the truncated count minus $O(\log n)$
and that count plus the displayed tail estimate.
Our choice of $J$ and sufficiently large $n$
therefore yield \eqref{eq:almost_all_R} for every
nonexceptional modulus. Finally, the divisibility
\eqref{eq:almost_all_full_power} follows from
Theorem~\ref{thm:main}.
\end{proof}

\newpage

\subsection{Expected Number of Twin Prime Pairs in \texorpdfstring{$C_n$}{Cn}}\label{sec:twins}
The prime number theorem in arithmetic progressions
makes no assertion about the existence of twin primes.
Nevertheless an actual twin-prime pair lying in the
prime interval of Lemma~\ref{lem:one} gives an explicit
divisor of $\sig(C_n)$.
For $n\ge2$, define
\[
 \mathcal T(n)=\{p:\ p,p+2\text{ prime},\ n+1<p<p+2\le2n\},
 \qquad T(n)=\#\mathcal T(n).
\]

\begin{proposition}[Twin-prime product]\label{prop:twins}
For every $n\ge2$,
\begin{equation}\label{eq:twins}
 \prod_{p\in\mathcal T(n)}(p+1)(p+3)\mid\sig(C_n),
 \qquad 24^{T(n)}\mid\sig(C_n).
\end{equation}
Consequently,
\[
 \vp{2}(\sig(C_n))\ge3T(n),
 \qquad \vp{3}(\sig(C_n))\ge T(n).
\]
\end{proposition}
\begin{proof}
Both members of each counted twin-prime pair belong
to $(n+1,2n]$, hence have exponent one in $C_n$.
Different counted pairs do not overlap: the only
possible overlapping twin-prime pairs are $(3,5)$
and $(5,7)$, and $(3,5)$ cannot be counted because
$n+1\ge3$. Multiplicativity and
Proposition~\ref{prop:product} give the first divisibility.
For any pair with $p>3$, write $p=6j-1$.
Then
\[
 (p+1)(p+3)=6j(6j+2)=12j(3j+1).
\]
Exactly one of $j$ and $3j+1$ is even, so this
product is divisible by $24$. Multiplying over
all pairs proves the second assertion; taking
prime-adic valuations proves the inequalities.
\end{proof}

For example, $C_9=4862=2\cdot11\cdot13\cdot17$ and
$\sig(C_9)=9072$. The pair $(11,13)$ contributes
$(11+1)(13+1)=168$, and $168\mid9072$.
The factor $24$ per pair is best possible at the
level of these pair contributions, since
$168=2^3\cdot3\cdot7$.

\begin{remark}
The conclusions of Proposition~\ref{prop:twins}
are unconditional, but make no assertion that
$T(n)>0$ for arbitrarily large $n$. Under the
Hardy--Littlewood twin-prime conjecture \cite{HardyLittlewood}
one would have $T(n)\sim 2C_2 n/(\log n)^2$,
where $C_2$ is the twin-prime constant. This
conditional prediction is not needed for any proof
in the paper; the resulting valuation bounds are
asymptotically weaker than Corollary~\ref{cor:Rasymp}.
\end{remark}

\section*{and the final remarks}
The initial divisor sums (with the indexing in
\eqref{eq:C}) are:
\begin{center}
\begin{tabular}{@{}r r r@{}}
\toprule
$n$ & $C_n$ & $\sig(C_n)$\\
\midrule
2 & 2 & 3\\
3 & 5 & 6\\
4 & 14 & 24\\
5 & 42 & 96\\
6 & 132 & 336\\
7 & 429 & 672\\
8 & 1430 & 3024\\
9 & 4862 & 9072\\
10 & 16796 & 35280\\
\bottomrule
\end{tabular}
\end{center}
Theorem~\ref{thm:main} supplies a rigorous
fixed-modulus divisibility result with a stronger
leading constant than the contribution from primes
in $(n+1,2n]$ alone.
Theorems~\ref{thm:siegel}--\ref{thm:nearoptimal} extend
the conclusion to growing moduli: uniformly for all
polylogarithmic moduli and for almost all moduli
up to nearly $\sqrt n$. In the almost-all result,
the exponent recovers the leading constant $2\log2$
from the full prime-interval decomposition. It does not specify an optimal
finite threshold for a given modulus. The
classification in Theorem~\ref{thm:classification}
and the exact gcd identities in
Proposition~\ref{prop:gcd} give the corresponding
universal assertions in exact form. The twin-prime contribution
is an explicit consequence of multiplicativity,
not a proof of any conjecture about the frequency
of twin primes.

The arguments rely on classical facts about primes
in arithmetic progressions and prime factors of
binomial coefficients. They clarify which divisor-sum
congruences follow directly from exponent-one prime
factors and where genuinely sharper information
would require an analysis of the higher prime-power
factors of $C_n$.

\bigskip
\begin{verse}\small
Ekmeğe, aşka ve ömre\\
Küfeleriyle hükmeden\\
Ciğerleri küçük, elleri büyük\\
Nefesleri yetmez avuçlarına\\
İlkokul çağında hepsi\\
Kenar çocukları\\
Kar altındadır\\
\emph{A. Arif}
\end{verse}
\end{document}